\documentclass[11pt]{article}

\usepackage{amscd,amsmath, amssymb, fancyhdr, graphicx, url}

\numberwithin{equation}{section}

\newcommand{\version}{version 1.0,\ \ Mar 4, 2015}

\makeatletter
\def\x@arrow{\DOTSB\Relbar}
\def\xlongrightarrowfill@{\arrowfill@\relbar\relbar\longrightarrow}
\newcommand{\xlongrightarrow}[2][]{%
        \ext@arrow 0099\xlongrightarrowfill@{#1}{#2}}
\makeatother

\def\eqref#1{(\ref{#1})}

\newcommand{\arrow}{{\:\longrightarrow\:}}
\newcommand{\Z}{{\Bbb Z}}
\def\C{{\Bbb C}}
\def\P{{\Bbb P}}

\newcommand{\R}{{\Bbb R}}
\newcommand{\Q}{{\Bbb Q}}

\def\1{\sqrt{-1}\:}

\newcommand{\cntrct}                
{\hspace{2pt}\raisebox{1pt}{\text{$\lrcorner$}}\hspace{2pt}}

\renewcommand{\tilde}{\widetilde}
\renewcommand{\bar}{\overline}
\renewcommand{\phi}{\varphi}
\renewcommand{\epsilon}{\varepsilon}

\newcommand{\Teich}{\operatorname{\sf Teich}}

\newcommand{\Comp}{\operatorname{\sf Comp}}
\newcommand{\Symp}{\operatorname{\sf Symp}}
\newcommand{\Hyp}{\operatorname{\sf Hyp}}

\newcommand{\Per}{\operatorname{\sf Per}}
\newcommand{\Perspace}{\operatorname{{\Bbb P}\sf er}}

\newcommand{\im}{\operatorname{im}}
\newcommand{\End}{\operatorname{End}}

\newcommand{\Gr}{\operatorname{Gr}}

\newcommand{\Id}{\operatorname{Id}}

\newcommand{\Pos}{\operatorname{Pos}}

\newcommand{\Diff}{\operatorname{\sf Diff}}

\renewcommand{\Re}{\operatorname{Re}}
\renewcommand{\Im}{\operatorname{Im}}

\newcounter{Mycounter}[section]
\newcounter{lemma}[section]
\newcounter{claim}[section]
\newcounter{sublemma}[section]
\newcounter{corollary}[section]
\renewcommand{\thecorollary}{{Corollary \thesection.\arabic{corollary}}}
\newcommand{\corollary}{%
    \setcounter{corollary}{\value{Mycounter}}
    \refstepcounter{corollary}
    \stepcounter{Mycounter}
    {\noindent \bf \thecorollary:\ }}

\newcounter{theorem}[section]
\renewcommand{\thetheorem}{{Theorem \thesection.\arabic{theorem}}}
\newcommand{\theorem}{%
    \setcounter{theorem}{\value{Mycounter}}
    \refstepcounter{theorem}
    \stepcounter{Mycounter}
    {\noindent \bf \thetheorem:\ }}

\newcounter{conjecture}[section]
\newcounter{proposition}[section]
\renewcommand{\theproposition}
      {{Proposition \thesection.\arabic{proposition}}}
\newcommand{\proposition}{%
    \setcounter{proposition}{\value{Mycounter}}
    \refstepcounter{proposition}
    \stepcounter{Mycounter}
    {\noindent \bf \theproposition:\ }}

\newcounter{definition}[section]
\renewcommand{\thedefinition}
      {{Definition~\thesection.\arabic{definition}}}
\newcommand{\definition}{%
    \setcounter{definition}{\value{Mycounter}}
    \refstepcounter{definition}
    \stepcounter{Mycounter}
    {\noindent \bf \thedefinition:\ }}

\newcounter{example}[section]
\newcounter{remark}[section]
\renewcommand{\theremark}{{Remark \thesection.\arabic{remark}}}
\newcommand{\remark}{%
    \setcounter{remark}{\value{Mycounter}}
    \refstepcounter{remark}
    \stepcounter{Mycounter}
    {\noindent \bf \theremark:\ }}

\newcounter{problem}[section]
\newcounter{question}[section]
\renewcommand{\thequestion}{{Question \thesection.\arabic{question}}}
\newcommand{\question}{%
    \setcounter{question}{\value{Mycounter}}
    \refstepcounter{question}
    \stepcounter{Mycounter}
    {\noindent \bf \thequestion:\ }}

\makeatletter

\@addtoreset{equation}{section} \@addtoreset{footnote}{section}
\makeatother

\def\blacksquare{\hbox{\vrule width 5pt height 5pt depth 0pt}}
\def\endproof{\blacksquare}

\begin{document}

\begin{center}
{\LARGE\bf
Teichm\"uller space for hyperk\"ahler\\[2mm] and symplectic structures\\[4mm]
}

Ekaterina Amerik\footnote{Partially supported by 
RScF grant, project 14-21-00053, 11.08.14.}, Misha Verbitsky\footnote{Partially supported by 
RScF grant, project 14-21-00053, 11.08.14.


}

\end{center}

{\small \hspace{0.15\linewidth}
\begin{minipage}[t]{0.7\linewidth}
{\bf Abstract} \\
Let $S$ be an infinite-dimensional manifold of all
symplectic, or hyperk\"ahler, structures on a compact
manifold $M$, and $\Diff_0$ the connected component of
its diffeomorphism group. The quotient $S/\Diff_0$ 
is called the Teichm\"uller space of symplectic
(or hyperk\"ahler) structures on $M$. 
MBM classes on a hyperk\"ahler
manifold $M$ are cohomology classes which can be represented
by a minimal rational curve on a deformation of 
$M$. We determine the Teichm\"uller
space of hyperk\"ahler structures on a hyperk\"ahler manifold,
identifying any of its connected components
with an open subset of the Grassmannian variety
$SO(b_2-3,3)/SO(3)\times SO(b_2-3)$ consisting of all
Beauville-Bogomolov positive 3-planes in $H^2(M, \R)$ which are not orthogonal 
to any of the MBM classes. This is used to determine
the Teichm\"uller space of symplectic structures
of 
K\"ahler type on a hyperk\"ahler manifold of maximal 
holonomy. We show that any connected component of 
this space is naturally identified with the space
of cohomology classes $v\in H^2(M,\R)$ with
$q(v,v)>0$, where $q$ is the 
Bogomolov-Beauville-Fujiki form on $H^2(M,\R)$.
\end{minipage}
}

\tableofcontents


\section{Introduction}


\subsection{Statement of the main result}

Denote by $\Gamma(\Lambda^2 M)$ the space of all
2-forms on a manifold $M$,
and let $\Symp\subset \Gamma(\Lambda^2 M)$ be the space of 
all symplectic 2-forms. We equip $\Gamma(\Lambda^2 M)$
with $C^\infty$-topology of uniform convergence 
on compacts with all derivatives. Then $\Gamma(\Lambda^2
M)$ is a Frechet vector space, and $\Symp$ a Frechet
manifold. We consider the group of diffeomorphisms, denoted
$\Diff$ or $\Diff (M)$, as a Fr\'echet Lie group, and denote
its connected component (known sometimes as the group
of isotopies) by $\Diff_0$. The quotient group
$\Gamma:=\Diff/\Diff_0$ is called {\bf the mapping class group}
of $M$.

Teichm\"uller space of symplectic structures on 
$M$ is defined as a quotient $\Teich_s:= \Symp/\Diff_0$.
It was studied in 
\cite{_PMH_Wilson:some_} and \cite{_Fricke_Habermann_}
together with its quotient $\Teich_s/\Gamma=\Symp/\Diff$,
known as {\bf the moduli space} of symplectic structures.

Notice that by Moore's theorem the action of $\Gamma$
on the space $\Teich_s(V)$ of symplectic forms
with fixed volume $V$ in $\Teich_s$ is ergodic e.g. for a compact torus of 
dimension $>2$ or a hyperk\"ahler manifold.\footnote{See 
\cite{_Entov_Verbitsky:packing_}, Theorem 9.2, where the argument is given in
all detail
.} Therefore, 
$\Gamma$ acts on $\Teich_s(V)$ with dense orbits, and the quotient
``space'' has a topology not much different from
the codiscrete one; in particular, it is not
``a manifold'', even in the most general sense
of this word. However, as 
shown in \cite{_Moser:volume_} 
(see also \cite{_Fricke_Habermann_}, Proposition 3.1), the space
$\Teich_s$ is a manifold, possibly non-Hausdorff,
and the symplectic period map
$\Per_s:\; \Teich_s \arrow H^2(M, \R)$, associating 
to $\omega \in \Teich_s$ its cohomology class,
is locally a diffeomorphism.

In the present paper we study the Teichm\"uller
space $\Teich_k$ of all symplectic structures of K\"ahler
type on a hyperk\"ahler manifold. 
Our main result is the
following theorem.

\hfill

\theorem\label{_teich_main_intro_Theorem_}
Let $M$ be a hyperk\"ahler manifold of maximal holonomy, and
$\Teich_k$ the space of all symplectic forms admitting
a compatible hyperk\"ahler structure (this is equivalent to being
K\"ahler in a certain complex structure, see \ref{_Calabi-Yau_Theorem_}).
Then the symplectic period map $\Per_s:\; \Teich_k \arrow H^2(M, \R)$
is an open embedding on each connected component, 
and its image is determined by quadratic 
inequality \[ \im\Per_s=\{v\in H^2(M,\R)\ \ |\ \ q(v,v)>0\},\]
where $q$ is the BBF form (\ref{_BBF_expli_Equation_}).

\hfill

{\bf Proof:} See \ref{th-teich-k}. \endproof

\hfill

The idea is to relate $\Teich_k$ to the Teichm\"uller space of 
hyperk\"ahler structures (\ref{def-teich-hk}) and to remark that
the latter has an explicit description in terms of {\bf MBM classes}
defined in \cite{_AV:MBM_}.

\subsection{The space of symplectic structures}

The study of the space of symplectic structures
was initiated by Moser in \cite{_Moser:volume_}.
Moser proved the following beautiful theorem, which
lies in foundation of symplectic topology.

\hfill

\theorem
Let $\omega_t$ be a continuous family of symplectic
structures on a compact manifold $M$. Assume that
the cohomology class of $\omega_t$ is constant
(that is, independent on $t$). Then all $\omega_t$
are related by diffeomorphisms.
\endproof

\hfill

It is not hard to see that this theorem implies
that the period map from the symplectic Teichm\"uller
space to cohomology is a local diffeomorphism (\cite{_Fricke_Habermann_}, Proposition 3.1). However, further
study of this space is very complicated, and in dimension
$>4$ almost nothing is known.

For a state of the art survey of the 
moduli of symplectic structures, please see 
\cite{_Salamon:symple_moduli_}.
For a particular interest to us is 
\cite[Example 3.3]{_Salamon:symple_moduli_},
due to D. McDuff (\cite{_MdD:examples_}).
It shows that the Teichm\"uller space of symplectic
structures on $S^2 \times S^2 \times T^2$ 
is non-Hausdorff.

The utility of our results for the general problem
of studying the symplectic structures is somewhat
restricted, because we consider only  symplectic structures
of K\"ahler type. It was conjectured, however, that
all symplectic structures on K3 are of K\"ahler
type (see e.g. \cite{_Donaldson:ellipt_}), hence this restriction
could theoretically be lifted. However, this conjecture
seems to be very hard. \ref{_teich_main_intro_Theorem_}
implies the following weaker form of this conjecture.

\hfill

\corollary
Let $M$ be a hyperk\"ahler manifold, 
and $\Teich_s$ a connected component of the Teichm\"uller
space of symplectic structures containing a
symplectic structure of hyperk\"ahler type.
Assume that $\omega\in \Teich_s$ is a symplectic
structure which is not of hyperk\"ahler type.
Then $\omega$ is a non-Hausdorff point in 
 $\Teich_s$, non-separable from a point
of hyperk\"ahler type.

\hfill

{\bf Proof:} Let $\Teich_k\subset \Teich_s$
be the set of all points of K\"ahler type in $\Teich_s$.
Then $\Per_s:\; \Teich_k \arrow H^2(M, \R)$
surjects to a connected component of all 
cohomology classes $\eta\in H^2(M)$ 
satisfying $\eta^{2n}\neq 0$, as follows from
\ref{_Fujiki_Theorem_}. After gluing all non-separable
points, the period map becomes an isomorphism,
since it is etale, and is an isomorphism
on $\Teich_k$. This implies that all 
points of $\Teich_s\backslash \Teich_k$
are non-separable from points in $\Teich_k$.
\endproof

\hfill

However, this does not imply that $\Teich_s=\Teich_k$,
because the symplectic Teichm\"uller space can be
non-Hausdorff, as follows from 
\cite[Example 3.3]{_Salamon:symple_moduli_}.

\hfill

\question
Let $\Teich_s$ be a  connected component of the Teichm\"uller
space of symplectic structures containing a
symplectic structure of hyperk\"ahler type.
Is it true that $\Teich_s$ is Hausdorff?

\hfill

Notice that the Teichm\"uller space of complex structures
on a hyperk\"ahler manifold is non-Hausdorff in all known examples.

\section{Hyperk\"ahler manifolds: basic results}

In this section, we recall the definitions and basic
properties of hyperk\"ahler manifolds and MBM classes.

\subsection{Hyperk\"ahler manifolds}

\definition
A {\bf hyperk\"ahler manifold}
is a compact, K\"ahler, holomorphically symplectic manifold.

\hfill

\definition
A hyperk\"ahler manifold $M$ is called
{\bf simple}, or {\bf maximal holonomy}, if $\pi_1(M)=0$ and $H^{2,0}(M)=\C$.

\hfill

This definition is motivated by the following theorem
of Bogomolov. 

\hfill

\theorem \label{_Bogo_deco_Theorem_}
(\cite{_Bogomolov:decompo_})
Any hyperk\"ahler manifold admits a finite covering
which is a product of a torus and several 
simple hyperk\"ahler manifolds.
\endproof

\hfill

The second cohomology $H^2(M,\Z)$ of a simple hyperk\"ahler manifold $M$ carries an integral quadratic form $q$, 
called {\bf the Bogomolov-Beauville-Fujiki form}. It was first
defined in \cite{_Bogomolov:defo_} and 
\cite{_Beauville_},
but it is easiest to describe it using the
Fujiki theorem, proved in \cite{_Fujiki:HK_}.

\hfill

\theorem\label{_Fujiki_Theorem_}
(Fujiki)
Let $M$ be a simple hyperk\"ahler manifold,
$\eta\in H^2(M)$, and $n=\frac 1 2 \dim M$. 
Then $\int_M \eta^{2n}=c q(\eta,\eta)^n$,
where 
$c>0$ is an integer. \endproof

\hfill

\remark 
Fujiki formula (\ref{_Fujiki_Theorem_}) 
determines the form $q$ uniquely up to a sign.
For odd $n$, the sign is unambiguously determined as well.
For even $n$, one needs the following explicit
formula, which is due to Bogomolov and Beauville.
\begin{equation}\label{_BBF_expli_Equation_}
\begin{aligned}
\tilde q(\eta,\eta'):=&
\int_M\eta\wedge\eta'\wedge\Omega^{n-1}\wedge\bar\Omega^{n-1}-\\&
-\frac{(n-1)}n\frac{\Big(\int_M\eta\wedge\Omega^{n-1}\wedge\bar\Omega^n
\Big)\cdot\Big(\int_M \eta'\wedge\Omega^n\wedge\bar\Omega^{n-1}\Big)}
{\int_M\Omega^n\wedge\bar\Omega^n},
\end{aligned}
\end{equation}
where $\Omega$ is the holomorphic symplectic form, and 
$\lambda>0$.

\subsection{MBM classes}



\definition
A cohomology class $\eta\in H^2(M, \R)$ is called 
{\bf positive} if $q(\eta,\eta)>0$, and 
{\bf negative}, if $q(\eta,\eta)<0$.

\hfill

\definition
Let $M$ be a hyperk\"ahler manifold. The 
{\bf monodromy group} of $M$ is a subgroup of $GL(H^2(M,\Z))$
generated by monodromy transforms for all Gauss-Manin local systems.
This group can also be characterized in terms of the
mapping class group action (\ref{_monodro_defi_Remark_}).

\hfill




The Beauville-Bogomolov-Fujiki form allows one to identify $H^2(M,\Q)$ and
$H_2(M,\Q)$. More precisely, it provides an embedding $H^2(M,\Z)\to
H_2(M,\Z)$ which is not an isomorphism (indeed $q$ is not necessarily
unimodular) but becomes an isomorphism after tensoring with $\Q$.
We thus can talk of classes of curves 
in $H^{1,1}(M, \Q)$, meaning
that the corresponding classes in $H_{1,1}(M, \Q)$ are classes of curves 
and shall do this in what follows.



Recall that {\bf a face} of a convex cone in a vector space $V$
is the intersection of its boundary and a hyperplane which 
has non-empty interior in the hyperplane.

If $f:(M,I)\dasharrow (M, I')$ is a birational isomorphism between
hyperk\"ahler manifolds, the induced map $f^*: H^2(M,I)\to H^2(M,I')$
is an isomorphism since $f$ is an isomorphism in codimension one.
By the K\"ahler cone of a birational model of $(M,I)$ as a part of $H^2(M,I)$
we mean the inverse image by $f$ of the K\"ahler cone of such an $(M,I')$.

\hfill

\definition\label{mbm}
A non-zero negative integral cohomology class
$z\in H^{1,1}(M,I)$ is called {\bf monodromy birationally minimal} (MBM)
if for some isometry $\gamma\in O(H^2(M,\Z))$ 
belonging to the monodromy group,
 $\gamma(z)^{\bot}\subset H^{1,1}(M,I)$ contains a face 
of the K\"ahler cone of one of birational
models $(M,I')$ of $(M,I)$.

\hfill

The following theorems summarize the main results about MBM classes from \cite{_AV:MBM_}.

\hfill

\theorem\label{defo-inv}(\cite{_AV:MBM_}, Corollary 5.13)
An MBM class $z\in H^{1,1}(M,I)$ is also MBM for any deformation of the complex structure $(M,I')$ where $z$ remains of
type $(1,1)$.

\hfill

\theorem\label{kahler-cone}(\cite{_AV:MBM_}, Theorem 6.2)
The K\"ahler cone of $(M,I)$ is 
a connected component of $\Pos(M,I)\backslash \cup_{z\in S} z^\bot$,
where $\Pos(M,I)$ is the positive cone of $(M,I)$ and $S$ is the set
of MBM classes on $(M,I)$.

\hfill

Because of the deformation-invariance property of MBM classes, it is natural to fix a connected component
$\Comp_0$ of the space of
complex structures of hyperk\"ahler type on $M$ and to call $z\in H^2(M,\Z)$ an MBM class (relative to $\Comp_0$)
when it is MBM in those complex structures where it is of type $(1,1)$. One moreover has the following description of
such classes.

\hfill

\corollary
A negative homology class $v\in H_2(M,\Z)$ is MBM if and only $\lambda v$
can be represented by an irreducible rational curve on $(M,J)$
for some $J\in \Comp_0$ and $\lambda\in \R^{\neq 0}$.

\hfill

{\bf Proof:} By deformation theory of hyperk\"ahler manifolds (see next section for some details, in 
particular \ref{zet-orthog}), there is a deformation $(M,J)$ 
of $(M,I)$ where only the multiples of $v$ survive 
as integral $(1,1)$-classes. By Theorem 5.15 of \cite{_AV:MBM_}, $v$ is MBM if and only if 
a multiple of $v$ is represented by a rational curve on $(M,J)$ (for reader's convenience
we recall that the reason behind this is that by the results of Huybrechts \cite{_Huybrechts:cone_} and Boucksom \cite{_Boucksom-cone_} 
the faces of the K\"ahler cone are given as orthogonals to classes of rational curves).
\endproof

\hfill









\section{Global Torelli theorem, hyperk\"ahler structures and monodromy group}


In this Section, we recall a number of 
results  about hyperk\"ahler manifolds,
used further on in this paper. For more
details and references, please see 
\cite{_Besse:Einst_Manifo_} and \cite{_V:Torelli_}.

\subsection{Hyperk\"ahler structures}
\label{_hk_twi_Subsection_}

\definition
Let $(M,g)$ be a Riemannian manifold, and $I,J,K$
endomorphisms of the tangent bundle $TM$ satisfying the
quaternionic relations
\[
I^2=J^2=K^2=IJK=-\Id_{TM}.
\]
The triple $(I,J,K)$ together with
the metric $g$ is called {\bf a hyperk\"ahler structure}
if $I, J$ and $K$ are integrable and K\"ahler with respect to $g$.

Consider the K\"ahler forms $\omega_I, \omega_J, \omega_K$
on $M$:
\[
\omega_I(\cdot, \cdot):= g(\cdot, I\cdot), \ \
\omega_J(\cdot, \cdot):= g(\cdot, J\cdot), \ \
\omega_K(\cdot, \cdot):= g(\cdot, K\cdot).
\]
An elementary linear-algebraic calculation implies
that the 2-form $\Omega:=\omega_J+\1\omega_K$ is of Hodge type $(2,0)$
on $(M,I)$. This form is clearly closed and
non-degenerate, hence it is a holomorphic symplectic form.

In algebraic geometry, the word ``hyperk\"ahler''
is essentially synonymous with ``holomorphically
symplectic'', due to the following theorem, which is
implied by Yau's solution of Calabi conjecture
(\cite{_Besse:Einst_Manifo_}, \cite{_Beauville_}).

\hfill

\theorem\label{_Calabi-Yau_Theorem_}
Let $M$ be a compact, K\"ahler, holomorphically
symplectic manifold, $\omega$ its K\"ahler form, $\dim_\C M =2n$.
Denote by $\Omega$ the holomorphic symplectic form on $M$.
Suppose that $\int_M \omega^{2n}=\int_M (\Re\Omega)^{2n}$.
Then there exists a unique hyperk\"ahler metric $g$ with the same
K\"ahler class as $\omega$, and a unique hyperk\"ahler structure
$(I,J,K,g)$, with $\omega_J = \Re\Omega$, $\omega_K = \Im\Omega$.
\endproof

\hfill

Further on, we shall speak of ``hyperk\"ahler manifolds''
meaning ``holomorphic symplectic manifolds of K\"ahler
type'', and ``hyperk\"ahler structures'' meaning
the quaternionic triples together with a metric. 

\hfill

Every hyperk\"ahler structure induces a whole 2-dimensional
sphere of complex structures on $M$, as follows. 
Consider a triple $a, b, c\in R$, $a^2 + b^2+ c^2=1$,
and let $L:= aI + bJ +cK$ be the corresponging quaternion. 
Quaternionic relations imply immediately that $L^2=-1$,
hence $L$ is an almost complex structure. 
Since $I, J, K$ are K\"ahler, they are parallel with respect
to the Levi-Civita connection. Therefore, $L$ is also parallel.
Any parallel complex structure is integrable, and K\"ahler.
We call such a complex structure $L= aI + bJ +cK$
{\bf a complex structure induced by the hyperk\"ahler structure}.
There is a 2-dimensional holomorphic family of 
induced complex structures, and the total space
of this family is called {\bf the twistor space}
of a hyperk\"ahler manifold, its base being {\bf the twistor line} in the Teichm\"uller space
$\Teich$ which we are going to define next.

\subsection{Global Torelli theorem and monodromy}

\definition
Let $M$ be a compact complex manifold, and 
$\Diff_0(M)$ a connected component of its diffeomorphism group
({\bf the group of isotopies}). Denote by $\Comp$
the space of complex structures on $M$, equipped with
its structure of a Fr\'echet manifold, and let
$\Teich:=\Comp/\Diff_0(M)$. We call 
it {\bf the Teichm\"uller space.}

\hfill

\remark
In many important cases, such as
for manifolds with trivial canonical class (\cite{_Catanese:moduli_}), 
$\Teich$ is a finite-dimensional
complex space; usually it is non-Hausdorff.

\hfill



\definition The {\bf mapping class group} is $\Gamma=\Diff(M)/\Diff_0(M)$. It
naturally acts on $\Teich$ (the quotient of $\Teich$ by $\Gamma$ may be
viewed as the ``moduli space'' for $M$, but in general it has very bad 
properties; see below).

\hfill

\remark
Let $M$ be a hyperk\"ahler manifold (as usually, we assume
$M$ to be simple). For any $J\in \Teich$,
$(M,J)$ is also a simple hyperk\"ahler manifold, 
because the Hodge numbers are constant in families.
Therefore, $H^{2,0}(M,J)$ is one-dimensional. 

\hfill

\definition
 Let 
\[ \Per:\; \Teich \arrow {\Bbb P}H^2(M, \C)
\]
map $J$ to the line $H^{2,0}(M,J)\in {\Bbb P}H^2(M, \C)$.
The map $\Per$ is 
called {\bf the period map}.

\hfill

\remark
The period map $\Per$ maps $\Teich$ into an open subset of a 
quadric, defined by
\[
\Perspace:= \{l\in {\Bbb P}H^2(M, \C)\ \ | \ \  q(l,l)=0, q(l, \bar l) >0\}.
\]
It is called {\bf the period domain} of $M$.
Indeed, any holomorphic symplectic form $l$
satisfies the relations $q(l,l)=0, q(l, \bar l) >0$,
as follows from \eqref{_BBF_expli_Equation_}.

\hfill

\proposition\label{_period_Grassmann_Proposition_}
The period domain $\Perspace$
is identified with the quotient
$SO(b_2-3,3)/SO(2) \times SO(b_2-3,1)$, which
is the Grassmannian of positive oriented 2-planes in $H^2(M,\R)$.

\hfill

{\bf Proof:} See for example \cite{_V:Torelli_}, section 2.4.





\hfill

\definition
Let $M$ be a topological space. We say that $x, y \in M$
are {\bf non-separable} (denoted by $x\sim y$)
if for any open sets $V\ni x, U\ni y$, $U \cap V\neq \emptyset$.

\hfill

\theorem\label{nonsep-birat} 
(Huybrechts; \cite{_Huybrechts:basic_}).
If two points $I,I'\in \Teich$ are non-separable, then  
there exists a bimeromorphism $(M,I)\dasharrow (M,I')$.
\endproof

\hfill


\hfill

\definition
The space $\Teich_b:= \Teich\!/\!\!\sim$ is called {\bf the
birational Teichm\"uller space} of $M$.

\hfill

\remark This terminology is slightly misleading since Teichm\"uller space is
not a moduli space and there are many Teichm\"uller points corresponding to
isomorphic manifolds (orbits of mapping class group action). In particular, there are 
non-separable points of the Teichm\"uller space which correspond
to biregular, not just birational, complex structures. Even for
K3 surfaces, the Teichm\"uller space is non-Hausdorff. 

\hfill

\theorem \label{_glo_Torelli_Theorem_}
(Global Torelli theorem; \cite{_V:Torelli_})
The period map 
$\Teich_b\stackrel \Per \arrow \Perspace$ is an isomorphism
on each connected component of $\Teich_b$.
\endproof

\hfill

\remark By a result of Huybrechts (\cite{_Huybrechts:finiteness_}), $\Teich$ has only finitely
many connected components. We shall denote by $\Teich_I$ the component 
containing the parameter point for the complex structure $I$, and by 
$\Gamma_I$
the subgroup of the mapping class group $\Gamma$ fixing this component.
Obviously $\Gamma_I$ is of finite index in $\Gamma$.

\hfill

\definition\label{_monodro_defi_Remark_}
The monodromy group of $(M,I)$ is the image of $\Gamma_I$ in the orthogonal
group $O(H^2(M, \Z), q)$.
As mentioned in the Introduction, it can also be described
as a subgroup of the group $O(H^2(M, \Z), q)$
generated by monodromy transform maps for 
Gauss-Manin local systems obtained from all
deformations of $(M,I)$ over a complex base
(\cite[Definition 7.1]{_V:Torelli_}). This is 
how this group was originally defined by Markman
(\cite{_Markman:constra_}, \cite{_Markman:survey_}).

\hfill

\remark 
A caution: usually ``the global Torelli theorem''
is understood as a theorem about Hodge structures.
For K3 surfaces, the Hodge structure on $H^2(M,\Z)$
determines the complex structure. 
For $\dim_\C M >2$, it is false.

\hfill

Finally, recall the following well-known
fact which shall be used in the sequel.

\hfill

\proposition\label{zet-orthog} Let $z\in H^2(M,\Z)$ 
be a cohomology class. The part of $\Teich$ corresponding
to the complex structures where $z$ is of type $(1,1)$ is the inverse image of 
$z^{\bot}\subset {\mathbb P}H^2(M,\C)$ (the orthogonal
being taken with respect to $q$). \endproof


\section{Teichm\"uller space of hyperk\"ahler structures}


\definition
Let $(M,I,J,K,g)$ and $(M,I',J',K',g')$ be two
hyperk\"ahler structures. We say that these structures
are {\bf equivalent} if the corresponding quaternionic algebras
in $\End(TM)$ coincide.

\hfill

\proposition\label{_hk_equiv_Proposition_}
Let $M$ be a hyperk\"ahler manifold,
and $(M,I,J,K,g)$ and $(M,I',J',K',g')$ be two
hyperk\"ahler structures of maximal holonomy. Then the following
conditions are equivalent.
\begin{description}
\item[(i)] $g$ is proportional to $g'$
\item[(ii)] $(M,I,J,K,g)$ is equivalent to
  $(M,I',J',K',g')$.
\end{description}
{\bf Proof of (i) $\Rightarrow$ (ii):} If $g$ is proportional to $g'$,
the corresponding Levi-Civita connections coincide.
The corresponding holonomy group is $Sp(n)$,
and its stabilizer in $\End(TM)$ is a quaternionic
algebra, generated by $I,J,K$ and also by $I',J',K'$.

{\bf Proof of (ii) $\Rightarrow$ (i):} Conversely, assume
that the algebras generated by $I,J,K$ and by $I',J',K'$
coincide. Recall that a manifold is called {\bf
hypercomplex} if it is equipped with a triple of 
complex structures $I,J,K$ satisfying quaternionic
relations. By Obata's theorem, there exists a
unique torsion-free connection preserving $I,J,K$
on any hypercomplex manifold, \cite{_Obata_}; this connection
is called {\bf the Obata connection}. Since the
Levi-Civita connection on $(M,I,J,K,g)$ satisfies
this condition, it coincides with the Obata
connection for $(M,I,J,K)$ and for $(M,I',J',K')$
(the latter is true because the corresponding 
hypercomplex manifolds are equivalent). However,
$g$ and $g'$ are invariant with respect to the
holonomy of Levi-Civita connection, which is equal to 
$Sp(n)$. The space of $Sp(n)$-invariant 
symmetric 2-forms is 1-dimensional (\cite{_Weyl:invariants_}), hence
$g$ and $g'$ are proportional.
\endproof

\hfill

\remark\label{equiv-h} It is easy to see that two hyperk\"ahler structures of maximal holonomy
$(M,I,J,K,g)$ and $(M,I',J',K',g')$ are equivalent if and only if there exists a unitary quaternion
$h$ such that $hIh^{-1}=I'$, $hJh^{-1}=J'$, $hKh^{-1}=K'$, and a positive constant $\lambda$ with 
$\lambda g=g'$.

Consider the infinite-dimensional space $\Hyp$ of all quaternionic
triples $I,J,K$ on $M$ which are induced by some
hyperk\"ahler structure, with the same $C^\infty$-topology
of convergence with all derivatives. The quotient \\ 
$\Hyp/SU(2)$ (which is probably better to write as 
$\Hyp/SO(3)$, since $-1$ acts trivially on the triples)  is naturally identified with the set
of equivalence classes of hyperk\"ahler structures, up to changing the metric $g$ by a constant.

\hfill

\remark\label{_metric_fixed_vol_hk_Remark_}
By \ref{_hk_equiv_Proposition_}, for manifolds
with maximal holonomy the quotient $\Hyp_m:=\Hyp/SU(2)$  is
also identified with the space of  all hyperk\"ahler
metrics of fixed volume, say, volume 1.

\hfill

\definition\label{def-teich-hk}
Define {\bf the Teichm\"uller space of 
hyperk\"ahler structures} as the quotient $\Hyp_m/\Diff_0$,
where $\Diff_0$ is the connected component of the group
of diffeomorphisms $\Diff$, and {\bf the moduli
of hyperk\"ahler structures} as $\Hyp_m/\Diff$. 

\hfill

\remark
For most geometric structures, the Teichm\"uller spaces 
and especially the moduli spaces are non-Hausdorff.
However, for manifolds of maximal holonomy
the moduli space of
hyperk\"ahler structures is Hausdorff. This is because
$\Hyp_m$ is the space of all hyperk\"ahler
metrics of fixed volume
(\ref{_metric_fixed_vol_hk_Remark_}).
However, there is a metric on the moduli 
space of all metrics, known as Gromov-Hausdorff
metric (\cite{_Gromov:Riemannian_}), 
and a metric space is necessarily Hausdorff.

\hfill

\remark
Let $(\omega_I, \omega_J, \omega_K)$
be a triple of classes obtained from a hyperk\"ahler structure
in a given component of $\Teich_h$. Then the 3-dimensional space
$\langle \omega_I, \omega_J, \omega_K\rangle$ is oriented. Indeed,
$\langle \omega_J, \omega_K\rangle$ determines a complex structure
$I$ by global Torelli theorem, and the sign of $\omega_I$
is determined by a component of $\Pos(M,I)$ containing
its K\"ahler cone.

\hfill

The main result of this section is the following
theorem.

\hfill

\definition\label{_Teich_h_Definition_}
Let $M$ be a hyperk\"ahler manifold of maximal holonomy,
and  $\Teich_h:= \Hyp_m/\Diff_0$ the Teichm\"uller
space of hyperk\"ahler structures. Consider the
space $\Perspace_h=Gr_{+++}(H^2(M,\R))$ of all positive 
oriented 3-dimensional subspaces
in $H^2(M,\R)$, naturally diffeomorphic to
$\Perspace_h\cong SO(b_2-3, 3)/SO(3)\times
SO(b_2-3)$. Let $\Per_h:\; \Teich_h \arrow\Perspace_h$ 
be the map associating to a hyperk\"ahler structure
$(M,I,J,K,g)$ the 3-dimensional space generated
by the three K\"ahler forms $\omega_I, \omega_J,\omega_K$.
This map called {\bf the period map for the 
Teichm\"uller space of hyperk\"ahler structures},
and $\Perspace_h$ {\bf the period space of  hyperk\"ahler
  structures}.

\hfill

\theorem\label{th-teich-hk}
Let $M$ be a hyperk\"ahler manifold of maximal holonomy,
and $\Per_h:\; \Teich_h \arrow\Perspace_h$ the period map
for the 
Teichm\"uller space of hyperk\"ahler structures.
Then the period map 
$\Per_h:\; \Teich_h \arrow\Perspace_h$ is an open embedding
for each connected component. Moreover, its
image is the set of all spaces $W\in \Perspace_h$
such that the orthogonal complement $W^\bot$
contains no MBM classes.

\hfill

{\bf Proof:} First we describe the image of the period map. 
Let $W\in \Perspace_h=Gr_{+++}H^2(M,\R)$ be the positive three-dimensional space 
corresponding to a point of $\Perspace_h$. Consider the set of positive planes in $W$: 
each plane $V$ in this set is an element of $Gr_{++}H^2(M,\R)=\Perspace$, thus a period point
of an irreducible holomorphic symplectic manifold $(M,I)$. If $W$ comes from an hyperk\"ahler
structure $(I,J,K)$, that is, $W=\langle \omega_I, \omega_J, \omega_K\rangle$, 
then the orthogonal to $V$ in $W$ 
is 
generated by the class of the K\"ahler form $\omega_I$, $V$ itself being generated by $Re(\Omega)$
and $Im(\Omega)$, where $\Omega$ is the holomorphic symplectic form on $(M, I)$ (see \ref{_Calabi-Yau_Theorem_}). If an integral class $z$ is orthogonal to $W$, it means that 
it is of type $(1,1)$ on $(M,I)$ for any $I$ corresponding to $V$ as above, by \ref{zet-orthog}. In particular
$z$ is orthogonal to $\omega_I$. Therefore $z$ cannot be MBM
 since the MBM classes are never orthogonal to K\"ahler classes.

Conversely, suppose that $W$ is not orthogonal to any MBM class. Then the same is true for
a sufficiently general plane $V\subset W$. This plane corresponds to the period point of an
irreducible holomorphic symplectic manifold $(M,I)$. Now take $v\in W$ orthogonal to $V$.
Up to a sign, this is an element of $\Pos(M,I)$. 

We remark that a cohomology class $\eta\in H^2(M)$ 
is of type $(1,1)$ with respect to $(M,I)$
if and only if $\eta$ is orthogonal to $V$.
Since $W$ is not orthogonal to any MBM class, and
$V\subset W$ is generic, no MBM class is of type
$(1,1)$ on $(M,I)$. This means that the 
K\"ahler cone of $(M,I)$ is equal to the positive
cone, that, is, up to a constant, $v$ is a K\"ahler class, and thus that there is a hyperk\"ahler
metric $g$ such that $v$ is proportional to $\omega_I$ and $W=\langle \omega_I, \omega_J, \omega_K \rangle$
(by \ref{_Calabi-Yau_Theorem_}). Therefore $W$ is in the image of $\Per_h$.

Finally, let us show that the period map is injective. As we have already mentioned, the
planes $V\subset W$ are period points of irreducible holomorphic symplectic manifolds,
parameterized by $S^2=\P^1$. Hyperk\"ahler structures correspond one-to-one to the twistor lines 
in the Teichm\"uller space $\Teich$. If two hyperk\"ahler structures $g_1, g_2$ give the same vector space
$W\in Gr_{+++}H^2(M,\R)$, the corresponding twistor lines $L_{g_1}, L_{g_2}$ have 
the same image in 
$\Perspace=Gr_{++}H^2(M,\R)$. However,
by the global Torelli theorem (\ref{_glo_Torelli_Theorem_}), 
this is only possible when
each point of $L_{g_1}$ is unseparable from some point of $L_{g_2}$. It is easy to see that
this never happens: indeed, as $W$ is not orthogonal to MBM classes, the K\"ahler cone
is equal to the positive cone at a general point of 
$L_{g_1}$ (as well as of $L_{g_2}$), and such points are Hausdorff points of $\Teich$
(see for example \cite{_Markman:survey_}).
\endproof


\section{Teichm\"uller space of symplectic structures}


Our goal in this section is to describe the Teichm\"uller space for symplectic structures of K\"ahler type
$\Teich_k$.

\hfill

\theorem\label{th-teich-k}
The symplectic period map $\Per_s: \Teich_k\to H^2(M,\R)$ 
is an embedding on each connected component, and its image
is the set of positive vectors in $H^2(M,\R)$.

\hfill

{\bf Proof. Step 1:} 
The period map is locally a diffeomorphism by Moser's theorem 
(\cite{_Moser:volume_}; see also \cite{_Fricke_Habermann_}, Proposition 3.1),
so we only have to show that it has connected 
fibers and describe the image. Describing the image
is easy: let $v$ be a positive vector in $H^2(M,\R)$, then we can always choose a positive 3-subspace
$W\subset H^2(M,\R)$ which contains $v$ and is not orthogonal to an MBM class. By the proof of 
\ref{th-teich-hk},
$W$ gives rise to an hyperk\"ahler structure such that $v$ is the class of $\omega_I$, therefore 
$v$ must be in the image of $\Per_s$; the inverse inclusion is clear since K\"ahler classes are 
positive. 

\hfill

{\bf Step 2:} To show that $\Per_s$ has connected fibers,
we consider the following diagram:
\begin{equation}\label{_Teich_CD_Diagram_}
\begin{CD}
 \widetilde{\Teich}_h @>P>> \Teich_k \\
@VV{\tilde \Per_h}V @V{\Per_s}VV\\
{\small \begin{array}{c} \{x,y,z\in H^2(M)|x^2=y^2=z^2>0, \\
       \text{$x,y,z$ is an oriented, orthogonal triple}
    \\ \text{$\langle x, y, z\rangle^\bot$ contains no MBM
      classes} \} 
\end{array}} @>F>> \{x \in H^2(M)|x^2>0\} \\
\end{CD}
\end{equation}
where  $\widetilde{\Teich}_h$ is 
the Teichm\"uller space for hyperk\"ahler
triples $(I,J,K)$ together with the metric $g$,
$P$ the forgetful map putting $(I,J,K,g)$
to the symplectic form $\omega_I$, and
$F$ the forgetful map putting $(x,y,z)$ to $x$.

\hfill

\remark 
The space $\widetilde{\Teich}_h$ can be considered as
an $SO(3)\times {\R}^+$-bundle over 
the Teichm\"uller space $\Teich_h$ 
introduced in \ref{_Teich_h_Definition_}.

\hfill

{\bf Step 3:} 
The horizontal arrows of
\eqref{_Teich_CD_Diagram_} are surjective (the upper one by Calabi-Yau theorem),
and $\tilde \Per_h$ is an isomorphism by
\ref{th-teich-hk}. Therefore, fibers of
$F$ are surjectively projected to the fibers 
of $\Per_s$. To prove that the latter are
connected, it suffices to show that
the fibers of the forgetful map $F$ 
are connected.

\hfill

{\bf Step 4:} 
The fibers of $F$ can be described as follows.
Let $x\in H^2(M,\R)$ be a positive vector,
and $x^\bot \subset H^2(M,\R)$ its orthogonal
complement. Then $F^{-1}(x)$ is the space of
oriented orthogonal pairs $y, z\in w^\bot$ such that
$x^2=y^2=z^2$, and $\langle x, y, z\rangle^\bot$ contains no 
MBM classes. This space is an $S^1$-fibration
over an open subset $X\subset \Gr_{++}(x^\bot)$
of the corresponding oriented Grassmannian
$\Gr_{++}(x^\bot)$ consisting of all 2-planes
not orthogonal to any of MBM classes in $x^\bot$.
Therefore, $X$ is a complement to a codimension
2 subset of those planes $W\in \Gr_{++}(x^\bot)$
which are orthogonal to some MBM class in $x^\bot$. 
A complement to a codimension
2 subset in a connected manifold is also connected. 
This implies that $X$, and hence $F^{-1}(x)$,
is connected, finishing the proof of 
\ref{th-teich-k}. \endproof

\hfill

{\bf Acknowledgements:}
We are grateful to M. Entov, A. Glutsuk 
and V. Shevchishin for interesting discussions
of the subject.

{

\noindent {\sc Ekaterina Amerik\\
{\sc Laboratory of Algebraic Geometry,\\
National Research University HSE,\\
Department of Mathematics, 7 Vavilova Str. Moscow, Russia,}\\
\tt  Ekaterina.Amerik@math.u-psud.fr}, also: \\
{\sc Universit\'e Paris-11,\\
Laboratoire de Math\'ematiques,\\
Campus d'Orsay, B\^atiment 425, 91405 Orsay, France}

\hfill

\noindent {\sc Misha Verbitsky\\
{\sc Laboratory of Algebraic Geometry,\\
National Research University HSE,\\
Department of Mathematics, 7 Vavilova Str. Moscow, Russia,}\\
\tt  verbit@mccme.ru}, also: \\
{\sc Kavli IPMU (WPI), the University of Tokyo}

}

\end{document}